\documentclass[a4paper,12pt]{amsart}
\usepackage{amsfonts}
\usepackage{amssymb}
\usepackage{ifthen}
\usepackage{graphicx}
\usepackage{color}
\nonstopmode \numberwithin{equation}{section}
\newtheorem{thm}{Theorem}%[section]
\newtheorem{lem}{Lemma}%[section]
\newtheorem{cor}{Corollary}%[section]
\newtheorem{prop}{Proposition}%[section]

\newtheorem{conj}{Conjecture}

\theoremstyle{definition}
\newtheorem{defn}{Definition}%[section]
\newtheorem{example}{Example}%[section]

\newtheorem{ques}{Question}
\newtheorem{rem}{Remark}

\newtheorem{rems}{Remarks}

\newcounter {own}
\def\theown {\thesection  .\arabic{own}}

\newenvironment{pf}[1][]{%
 \vskip 3mm
 \noindent
 \ifthenelse{\equal{#1}{}}%
  {{\slshape Proof. }}%
  {{\slshape #1.} }%
 }%
{\qed\bigskip}

\newcounter{alphabet}
\newcounter{tmp}

\newcommand{\IN}{{\mathbb N}}
\newcommand{\IC}{{\mathbb C}}

\newcommand{\D}{{\mathbb D}}

\newcommand{\dist}{{\operatorname{dist}}}

\def\be{\begin{equation}}
\def\ee{\end{equation}}

\newcommand{\bee}{\begin{enumerate}}
\newcommand{\eee}{\end{enumerate}}

\newcommand{\blem}{\begin{lem}}
\newcommand{\elem}{\end{lem}}
\newcommand{\bthm}{\begin{thm}}
\newcommand{\ethm}{\end{thm}}
\newcommand{\bcor}{\begin{cor}}
\newcommand{\ecor}{\end{cor}}
\newcommand{\beg}{\begin{example}}
\newcommand{\eeg}{\end{example}}
\newcommand{\begs}{\begin{examples}}
\newcommand{\eegs}{\end{examples}}
\newcommand{\bdefn}{\begin{defn}}
\newcommand{\edefn}{\end{defn}}
\newcommand{\bprob}{\begin{prob}}
\newcommand{\eprob}{\end{prob}}
\newcommand{\bei}{\begin{itemize}}
\newcommand{\eei}{\end{itemize}}
\newcommand{\bqn}{\begin{ques}}
\newcommand{\eqn}{\end{ques}}
\newcommand{\bcon}{\begin{conj}}
\newcommand{\econ}{\end{conj}}
\newcommand{\bcons}{\begin{conjs}}
\newcommand{\econs}{\end{conjs}}
\newcommand{\bprop}{\begin{prop}}
\newcommand{\eprop}{\end{prop}}
\newcommand{\brem}{\begin{rem}}
\newcommand{\erem}{\end{rem}}
\newcommand{\brems}{\begin{rems}}
\newcommand{\erems}{\end{rems}}
\newcommand{\bo}{\begin{obser}}
\newcommand{\eo}{\end{obser}}
\newcommand{\bos}{\begin{obsers}}
\newcommand{\eos}{\end{obsers}}
\newcommand{\bpf}{\begin{pf}}
\newcommand{\epf}{\end{pf}}
\newcommand{\ba}{\begin{array}}
\newcommand{\ea}{\end{array}}
\newcommand{\beq}{\begin{eqnarray}}
\newcommand{\beqq}{\begin{eqnarray*}}
\newcommand{\eeq}{\end{eqnarray}}
\newcommand{\eeqq}{\end{eqnarray*}}

\newcommand{\ra}{\rightarrow}

\newcommand{\ds}{\displaystyle}

\newcounter{minutes}
\divide\time by 60
\newcounter{hours}
\multiply\time by 60 \addtocounter{minutes}{-\time}
\begin{document}
\bibliographystyle{amsplain}
\title[Discrete Analogue of Generalized Hardy Spaces and Multiplication Operators on Homogenous Trees]
{Discrete Analogue of Generalized Hardy Spaces and Multiplication Operators on Homogenous Trees}

%%%=========================================================================
%%\thanks{%$^\dagger$
%File:~\jobname .tex,
%          printed: \number\day-\number\month-\number\year,
%          \thehours.\ifnum\theminutes<10{0}\fi\theminutes}
%%%=========================================================================
%
%\author{S. Ponnusamy $^\dagger $}

\author{Perumal Muthukumar}
\address{P. Muthukumar,
Indian Statistical Institute (ISI), Chennai Centre,
SETS (Society for Electronic Transactions and Security),
MGR Knowledge City, CIT Campus, Taramani,
Chennai 600 113, India. }
\email{pmuthumaths@gmail.com}

\author{Saminathan Ponnusamy  %$^\dagger $
}
\address{S. Ponnusamy,
Indian Statistical Institute (ISI), Chennai Centre,
SETS (Society for Electronic Transactions and Security),
MGR Knowledge City, CIT Campus, Taramani,
Chennai 600 113, India.
}
\email{samy@isichennai.res.in, samy@iitm.ac.in}

\subjclass[2000]{Primary: 47B38, 05C05; Secondary: 46B50,46B26}
\keywords{Rooted homogeneous tree, multiplication operators, Lipschitz space, generalized Hardy spaces\\
%$^\dagger$ {\tt Corresponding author}\\
%This first author is currently on leave from the Department of Mathematics, Indian Institute of Technology Madras, Chennai-600 036, India.
}

%\date{\today
%File: Muthu.tex}

\begin{abstract}
In this article, we define discrete analogue of generalized Hardy spaces and its separable subspace on a homogenous rooted tree and
study some of its properties such as completeness, inclusion relations with other spaces, separability, growth estimate for functions
in these spaces and their consequences. Equivalent conditions for multiplication operators to be bounded and compact are also obtained.
Furthermore, we discuss about point spectrum, approximate point spectrum and spectrum of multiplication operators and discuss
when a multiplication operator is an isometry.
\end{abstract}
\thanks{
%%=========================================================================
File:~\jobname .tex,
          printed: \number\day-\number\month-\number\year,
          \thehours.\ifnum\theminutes<10{0}\fi\theminutes
%%=========================================================================
}
\maketitle
\pagestyle{myheadings}
\markboth{P. Muthukumar and S. Ponnusamy}{Generalized Hardy Spaces and Multiplication Operators on Homogenous Trees}

\section{Introduction}
The theory of function spaces defined on the unit disk ${\mathbb D}=\{z \in {\mathbb C}:\, |z|<1\}$ is particulary a well developed subject.
The recent book by Pavlovi\'c \cite{pavlovic:Book} and the book of Zhu \cite{Zhu:Book} provide us with a solid foundation
in studying various function spaces on the unit disk. One can also refer \cite{Duren:Hpspace} for Hardy spaces $(H^p)$,
\cite{Zhu:Bergmanspace} for Bergman spaces $(A^p)$, \cite{Book:Dirichlet} for Dirichlet spaces $(\mathcal{D}_p)$ and \cite{Bloch:functions} for Bloch space
$(\mathcal B)$.

In recent years, there has been a considerable interest in the study of function spaces on discrete set such as tree (more generally on graphs).
For example, Lipschitz space of a tree (discrete analogue of Bloch space) \cite{Colonna-MO-1}, weighted Lipschitz space of a tree \cite{Colonna-MO-2}, iterated logarithmic Lipschitz space of a tree \cite{Colonna-MO-3} and $H^p$ spaces on trees \cite{Hpspace:tree} are some in this line of investigation.
In \cite{Hpspace:tree} the $H^p$ spaces on trees are defined by means of certain maximal or square function operators associated with a nearest neighbour transition operator which is very regular and this study was further developed in \cite{Hpspace:green}.

In the study of operators on function spaces, multiplication and composition operators arise naturally and play an important role. Moreover, study of composition operators on various function spaces on the unit disk is
better developed than the multiplication operators, and the literature in these topics are exhaustive. See
for example the survey articles \cite{survey:Hp-multipliers, CO-MOsurvey, Dragan:survey} on multiplication operators on various function spaces of the
unit disk.

In the case of operator theory on discrete function spaces, average operator (Laplacian operator) is studied more than other operators.
Colonna and others studied multiplication operators on Lipschitz space, weighted Lipschitz space, and iterated logarithmic Lipschitz space of
a tree in \cite{Colonna-MO-1}, \cite{Colonna-MO-2} and \cite{Colonna-MO-3}, respectively. In \cite{Colonna-MO-5, Colonna-MO-4},
the authors discussed about multiplication operators between Lipschitz type spaces and the space of bounded functions on a tree.
Recently, Allen et al. \cite{Colonna-CO} have also studied composition operators on the Lipschitz space of a tree.

In this article, we define discrete analogue of generalized Hardy spaces and study its important properties.
Also, we study multiplication operators on generalized Hardy spaces on homogeneous trees
(i.e. every vertex has same number of neighbours).

\section{Preliminaries}
A \textit{graph} $G$ is a pair $G=(V,E)$ of sets satisfying $E\subseteq V \times V$. The elements of $V$ and $E$ are called
vertices and edges of the graph $G$, respectively. Two vertices $x,y \in V$ (with the abuse of language, one can write
as $x,y \in G$) are said to be \textit{neighbours} or adjacent (denoted by $x\sim y$) if there is an edge connecting them.
A \textit{regular} (homogeneous) graph is a graph in which every vertex has the same number of neighbours.
If every vertex has $k$ neighbours, then the graph is said to be \textit{$k-$regular} ($k-$homogeneous) graph.
A \textit{path} is part of a graph with finite or infinite sequence of distinct vertices $[v_0,v_1,v_2,\ldots]$ such
that $v_n\sim v_{n+1}$. If $P=[v_0-v_1-v_2-\cdots -v_{n}]$ is a path then the graph $C=[v_0-v_1-v_2-\cdots-v_{n}-v_0]$
(path $P$ with an additional edge $v_{n}v_{0}$) is called a \textit{cycle}. A non-empty graph $G$ is called \textit{connected}
if any two of its vertices are linked by a path in $G$. A connected and locally finite (every vertex has finite number of neighbours)
graph without cycles, is called a \textit{tree}. A \textit{rooted tree} is a tree in which a special vertex (called root) is
singled out. The \textit{distance} between any two vertex of a tree is the number of edges in the unique path connecting them.
If $G$ is a rooted tree with root $\textsl{o}$, then $|v|$ denotes the distance between $\textsl{o}$ and $v$. Further
the \textit{parent} (denoted by $v^-$) of a vertex $v$ is the unique vertex $w\in G$ such that $w\sim v$ and $|w|=|v|-1$.
For basic issues concerning graph theory, one can refer to standard texts such as \cite{Book:graph}.

Let $T$ be a rooted tree. By a function defined on a graph, we mean a function defined on its vertices. The Lipschitz space and
the weighted Lipschitz space of $T$ are denoted by $\mathcal{L}$ and $\mathcal{L}_w$, respectively. These are defined as follows:
$$ \mathcal{L} =\left \{f\colon T\to\mathbb{C}: \sup\limits_{v\in T,\,v\neq \textsl{o}} |f(v)-f(v^-)| <\infty \right \}
$$
and
$$ \mathcal{L}_w=\left \{f\colon T\to\mathbb{C}: \sup\limits_{v\in T,v\neq \textsl{o}} |v|\,|f(v)-f(v^-)| <\infty \right \},
$$
respectively. Throughout the discussion, a homogeneous rooted tree with infinite vertices is denoted by $T$,
$\mathbb{N}=\{1,2,\ldots \}$ and $\mathbb{N}_0=\mathbb{N}\cup \{0\}$.

For $p\in(0,\infty]$, the {\it generalized Hardy space
$H^{p}_{g}(\mathbb{D})$} consists of all those measurable functions
$f:\mathbb{D}\rightarrow\mathbb{C}$ such that $M_{p}(r,f)$ exists for all $r\in [0,1)$ and $\|f\|_{p}<\infty$,
where
$$\|f\|_{p}=
\begin{cases}
\displaystyle\sup_{0\leq r<1}M_{p}(r,f)
& \mbox{if } p\in(0,\infty),\\
\displaystyle\sup_{z\in\mathbb{D}}|f(z)| &\mbox{if } p=\infty,
\end{cases}
$$
and
$$
M_{p}^{p}(r,f)=\frac{1}{2\pi}\int_{0}^{2\pi}|f(re^{i\theta})|^{p}\,d\theta.
$$
The classical Hardy space $H^{p}$ is a subspace of $ H^{p}_{g}(\mathbb{D})$ consists of only analytic functions. See
\cite{Samy} for recent investigation on $ H^{p}_{g}(\mathbb{D})$ and some related function spaces.

For our investigation this definition has an analog in the following form.

\bdefn\label{def1}
Let $T$ be a $q+1$ homogeneous tree rooted at $\textsl{o}$. For every $n\in \mathbb{N}$, we introduce
$$%\be\label{eq2}
M_{p}(n,f):=
\left\{
\begin{array}{ll}
\ds \left (\frac{1}{(q+1)q^{n-1}}\sum\limits_{|v|=n}|f(v)|^{p} \right )^{\frac{1}{p}} & \mbox{if } p\in(0,\infty) \\
\max\limits_{|v|=n } |f(v)| & \mbox{if } p=\infty ,
\end{array}
\right.
$$
$M_{p}(0,f):= |f(\textsl{o})|$ and
\be\label{eq1}
\|f\|_{p}:= \sup\limits_{n\in \mathbb{N}_{0}} M_{p}(n,f).
\ee

The discrete analogue of the generalized Hardy space, denoted by $\mathbb{T}_{q,p}$, is then defined by
$$\mathbb{T}_{q,p}:=\{f\colon T \to\mathbb{C}~\mbox{ such that }~ \|f\|_{p}<\infty\}.
$$
Similarly, the discrete analogue of the generalized little Hardy space, denoted by $\mathbb{T}_{q,p,0}$, is defined by
$$\mathbb{T}_{q,p,0}:=\{f\in\mathbb{T}_{q,p} :\, \lim\limits_{n\rightarrow\infty} M_{p}(n,f)=0 \}
$$
for every $p\in(0,\infty]$. For the sake of simplicity, we shall write $\mathbb{T}_{q,p}$ and $\mathbb{T}_{q,p,0}$ as
$\mathbb{T}_{p}$ and $\mathbb{T}_{p,0}$, respectively.
Unless otherwise stated explicitly, throughout $\|.\|_p$ is defined as above.
\edefn

The paper is organized as follows.  In several sub sections of Section \ref{Main results},  we investigate important properties
such as completeness, inclusion relations and separability of the spaces
$\mathbb{T}_{p}$ and $\mathbb{T}_{p,0}$. First, we prove that $\mathbb{T}_{p}$ and $\mathbb{T}_{p,0}$ are Banach spaces
(see Theorems \ref{thm:banachp} and \ref{thm:banachp0}).
As with the $H^p$ spaces, similar inclusion relations are derived for the discrete cases (see Theorem \ref{cor1}).
Also, we show that $\mathbb{T}_{p}$ is not a separable space whereas $\mathbb{T}_{p,0}$ is separable (see Theorems \ref{thm:sepp0} and  \ref{thm:separable}).
Finally, in Section \ref{Main results}, we compare the convergence in norm  and uniform convergence on compact sets of $T$.
In Section \ref{Multiplication Operators},  we study multiplication operators on generalized Hardy spaces on homogeneous trees. In particular,
we give characterization of bounded and compact multiplication operators
on $\mathbb{T}_{p}$  (resp. on $\mathbb{T}_{p,0}$), see Theorems \ref{thm:boundedness} and \ref{thm:Compactness}. Also we determine the point spectrum,
the approximate point spectrum and the spectrum for the multiplication operator $M_\psi$ and  determine an upper bound for the essential
norm of $M_\psi$. Finally, we provide necessary and sufficient condition for a multiplication operator
to be an isometry.

\section{Topological Properties of $\mathbb{T}_{p}$ and $\mathbb{T}_{p,0}$} \label{Main results}

\subsection{Completeness}

\bthm\label{thm:banachp}
For $1\leq p\leq\infty$, $\|.\|_{p}$ induces a Banach space structure on the space $\mathbb{T}_{p}$.
\ethm
\bpf
First we begin with the case $p=\infty$. In this case, \eqref{eq1} reduces to $\|f\|_{p}=\sup\limits_{v\in T}|f(v)|$ and thus, the space $\mathbb{T}_{\infty}$
coincides with the set of all bounded functions on $T$ with sup-norm which is known to be a Banach space.

Next, we consider the case $1\leq p <\infty$. We have the following.
\begin{enumerate}
\item[(i)] If $f\equiv 0$, then $\|f\|_{p}=0$. Conversely, if $\|f\|_{p}=0$ then $M_{p}(n,f)=0$ for all $n\in \mathbb{N}_{0}$
showing that $\sum\limits_{|v|=n}|f(v)|^{p}=0$ for all $n\in \mathbb{N}_{0}$
and thus, $f\equiv 0$.

\item[(ii)] For each $n\in \mathbb{N}_0$ and $\alpha\in\mathbb{C}$, it is easy to see by definition that $M_{p}(n,\alpha f)= |\alpha|M_{p}(n,f)$ and thus,
$\|\alpha f\|_{p}= |\alpha|\|f\|_{p}$.

\item[(iii)] For each $n\in \mathbb{N}$ and $f,g \in \mathbb{T}_{p}$, one has (since $p\geq 1$)
\beqq
M_{p}(n,f+g) & =& \left \{\frac{1}{(q+1)q^{n-1}}\sum\limits_{|v|=n}(|f(v)|+|g(v)|)^{p} \right \}^{\frac{1}{p}} \\
&\leq & \left (\frac{1}{(q+1)q^{n-1}}\sum\limits_{|v|=n}|f(v)|^{p} \right )^{\frac{1}{p}}
+ \left (\frac{1}{(q+1)q^{n-1}}\sum\limits_{|v|=n}|g(v)|^{p} \right )^{\frac{1}{p}}\\
&=& M_{p}(n,f)+ M_{p}(n,g).
\eeqq
The last inequality trivially holds for $n=0$ and thus, $\|f+g\|_{p}\leq\|f\|_{p} + \|g\|_{p}$.
\end{enumerate}
Hence $(\mathbb{T}_{p}, \|.\|_{p})$ is a normed linear space. In order to prove that $\mathbb{T}_{p}$ is a Banach space,
we begin with a Cauchy sequence $\{f_k\}$ in $\mathbb{T}_{p}$. Then $\{f_k(v)\}$ is a Cauchy sequence in $\mathbb{C}$ for
every $v\in T$ and thus, $\{f_k\}$ converges pointwise to a function $f$. Now, for a given $\epsilon > 0$,
there exists an $N\in\mathbb{N}$ such that $M_{p}(n,f_k-f_l)<\epsilon$ for all $k,l\geq N$ and $n\in \mathbb{N}_{0}$.
Letting $l\rightarrow \infty$, we get $M_{p}(n,f_k-f)\leq\epsilon$ for all $k\geq N$ and $n\in \mathbb{N}_{0}$.
Hence $\|f_k-f\|_{p}\leq \epsilon$ for all $k\geq N$, which gives that $ f_k \rightarrow f.$ The triangle inequality
$\|f\|_{p} \leq \|f-f_N\|_{p} + \|f_N\|_{p}$ gives that $f\in \mathbb{T}_{p}$. This completes the proof of the theorem.
\epf

\brem
 For $0<p<1, \mathbb{T}_{p}$ is obviously a complete metric space.
\erem

A function $f\colon T \to\mathbb{C}$ is said to be \textit{radial constant function} if $f(v)=f(w)$ whenever $|v|=|w|$.

\brem
Since the integral means $M_{p}(r,f)$ of an analytic function $f$ defined on $\D$ is an increasing function of $r$, we have
$\|f\|_{p} = \lim\limits_{r\rightarrow 1^-}M_{p}(r,f)$. Thus, the little Hardy space $H_0^p$, defined by
$$H_0^p:= \{f\in H^p: M_{p}(r,f) \rightarrow 0 ~\mbox{ as }~ r\rightarrow 1^-\},$$
consists of only a single element, namely, the zero function. But this is not the case in the generalized Hardy space $H^{p}_{g}(\mathbb{D})$
of measurable functions, This is because the maximum modulus principle is not valid
for a general element in $H^{p}_{g}(\mathbb{D})$. Consequently, the generalized little Hardy space $H_{0,g}^p$ is
non-trivial (i.e. not a zero subspace), where
$$H_{0,g}^p:= \{f\in H^{p}_{g}(\mathbb{D}): M_{p}(r,f) \rightarrow 0 ~\mbox{ as }~ r\rightarrow 1^-\}.
$$
For example,  if
$$ f_{\alpha}(z)= \left \{\begin{array}{cl}
1 & \mbox{ for $|z|\leq \alpha$  }\\
0 & \mbox{ for $\alpha<|z|<1,$}
\end{array}
\right.
$$
then $f\in H^{p}_{0,g}$ for each $0\leq \alpha <1.$

In the discrete case, $\mathbb{T}_{p,0}$ is non-trivial. In fact, the set of all radial constant functions
in $\mathbb{T}_{p,0}$ is isometrically isomorphic to the sequence space $c_0$ (set of all sequences that converges to zero).
\erem

Here are few questions that arise naturally.

\bqn
Is it possible to define an inner product so that $\mathbb{T}_{2}$ becomes a Hilbert space? As with the $l^p$ and the $H^p$ spaces,
whether $\mathbb{T}_{p}$ is not isomorphic to $\mathbb{T}_{q}$ when $p\neq q$? What can be said about the dual of $\mathbb{T}_{p}$?
\eqn

These questions are open for the moment.

\bthm\label{thm:banachp0}
For $1\leq p\leq\infty$, $\|.\|_{p}$ induces a Banach space structure on $\mathbb{T}_{p,0}$.
\ethm
\bpf
For $n\in \mathbb{N}_{0}$ and $f,g \in \mathbb{T}_{p}$, we easily have
$$M_{p}(n,\alpha f)= |\alpha|M_{p}(n,f) ~\mbox{ and }~ M_{p}(n,f+g)\leq M_{p}(n,f)+ M_{p}(n,g)
$$
so that $\mathbb{T}_{p,0}$ is a subspace of $\mathbb{T}_{p}$. Suppose $\{f_k\}$ is a Cauchy sequence in $\mathbb{T}_{p,0}$.
Since $\mathbb{T}_{p}$ is a Banach space, $\{f_k\}$ converges to some function $f\in \mathbb{T}_{p}$. Next, we need to prove
that $f\in \mathbb{T}_{p,0}$, i.e., $M_{p}(n,f)\rightarrow 0$ as $n\rightarrow\infty$.
To do this, let $\epsilon > 0$ be given. Then there exists a $k\in \mathbb{N}$ such that $\|f_k-f\|_{p}< \epsilon/2$.
Since $f_k\in \mathbb{T}_{p,0}$, we can choose $N\in \mathbb{N}$ so that $M_p(n,f_k)<\epsilon/2$ for all $n\geq N$.
From the inequality $M_{p}(n,f)\leq M_{p}(n,f-f_k)+ M_{p}(n,f_k)$,
it follows that $M_{p}(n,f)<\epsilon$ for all $n\geq N$. Thus, $f\in \mathbb{T}_{p,0}$ which completes the proof.
\epf

\subsection{Inclusion Relations}

\blem\label{lem:inclusion}
For $0<r<s\leq\infty$ and for every complex valued function $f$ on $T$, $M_r(n,f)\leq M_s(n,f)$ holds for all $n\in \mathbb{N}_{0}$.
\elem
\bpf
The result for $s=\infty$ follows from the definition of $M_r(n,f)$ and thus, it suffices to prove the lemma for the case $0<r<s<\infty$.
Again by Definition \ref{def1}, we see that
$$M_r(0,f)=|f(\textsl{o})|=M_s(0,f).
$$
For $n\in\IN$, we have $(q+1)q^{n-1}$ vertices with $|v|=n$. Recall that on the Euclidean space $\IC^N$, the following norm equivalence
is well-known for $0<r<s$:
\begin{equation}\label{norm equivalence}
\|x\|_s\leq \|x\|_r\leq N^{\frac{1}{r}-\frac{1}{s}}\|x\|_s,
\end{equation}
where $p$-norm $\|.\|_p$ on $\IC^N$ is given by $\|x\|_p^p=\sum_{k=1}^N|x_k|^p$. The second inequality in \eqref{norm equivalence},
is an easy consequence of H\"{o}lder's inequality for finite sum. We may now use this with $N=(q+1)q^{n-1}$. As a consequence of it,
we have
$$\left (\sum\limits_{|v|=n}|f(v)|^{r}\right )^{\frac{1}{r}}
\leq \left \{(q+1)q^{n-1}\right \}^{\frac{1}{r}-\frac{1}{s} }\left (\sum\limits_{|v|=n}|f(v)|^{s}\right )^{\frac{1}{s}}
$$
which may be rewritten as
$$\left (\frac{1}{(q+1)q^{n-1}}\sum\limits_{|v|=n}|f(v)|^{r} \right )^{\frac{1}{r}} \leq
\left (\frac{1}{(q+1)q^{n-1}}\sum\limits_{|v|=n}|f(v)|^{s} \right )^{\frac{1}{s}}.
$$
This shows that $M_r(n,f)\leq M_s(n,f)$ for all $n\in \mathbb{N}_{0}$. The proof is complete.
\epf

As an immediate consequence of Lemma \ref{lem:inclusion}, one has the following.

\bthm\label{cor1}
For $0<r<s\leq\infty$, we have $\mathbb{T}_{s} \subset \mathbb{T}_{r}$
and $\mathbb{T}_{s,0}\subset \mathbb{T}_{r,0}$.
\ethm

We now show by an example that the inclusions in Theorem \ref{cor1} is proper.
Let $0<r<p<s\leq \infty$. Choose a sequence of vertices $\{v_n\}$ such that $|v_n|=n$ for all $n\in \mathbb{N}$. Consider the function $f$
defined by
$$f(v)= \left \{\begin{array}{cl}
\ds \{(q+1)q^{n-1}\}^{\frac{1}{p}} & \mbox{ if $v=v_n$ for some $n \in \mathbb{N}$ }\\
0 &\mbox{ elsewhere.}
\end{array}
\right .
$$
Then, for $n\in \IN$, one has
$$M_r(n,f)=\left \{(q+1)q^{n-1}\right \}^{ \frac{1}{p}-\frac{1}{r} }
$$
so that $M_r(n,f)\ra0$ as $n\ra\infty$, since $p>r$. Also, we have
$$M_s(n,f)= \left \{\begin{array}{ll}
\ds \{(q+1)q^{n-1}\}^{\frac{1}{p}} & \mbox{ if $s=\infty$}\\
\ds \{(q+1)q^{n-1}\}^{\frac{1}{p}-\frac{1}{s}} &\mbox{ if $s<\infty$ }
\end{array}
\right .
$$
and in either case, we find that $M_s(n,f)\ra\infty$ as $n\ra\infty$. This example shows that $\mathbb{T}_{s,0}$ is a proper subspace of $\mathbb{T}_{r,0}$.
From this example, it can be also seen that $\mathbb{T}_{s}$ is a proper subspace of $\mathbb{T}_{r}$.

\brems
\begin{enumerate}
\item
The unbounded function $f(v)=|v|$ belongs to $\mathcal{L}$ but is not in $\mathbb{T}_{p}$ for any $0<p<\infty$. On the other hand,
let us now fix an infinite path $\textsl{o}-v_1-v_2 \cdots $ with $|v_k|=k$. Define $g(v)=\{(q+1)q^{k-1}\}^{\frac{1}{p}}$
if $v=v_k$ and $0$ otherwise. It is then easy to check that $g$ belongs to $\mathbb{T}_{p}$ for all $0<p<\infty$ but is not in $\mathcal{L}$.
Thus, $\mathcal{L}$ is not comparable with $\mathbb{T}_{p}$ for all $0<p<\infty$.
\item
Consider the radial constant function $h$ defined by $h(\textsl{o})=0$ and
$$h(v)=\sum\limits_{k=1}^{|v|} \frac{1}{k} ~\mbox{ if $|v|\geq 1$}.
$$
By simple calculations, we find that $h$ belongs to $\mathcal{L}_w$ but is not in $\mathbb{T}_{p}$ for any $0<p<\infty$.
For the other direction, we fix an infinite path $\textsl{o}-v_1-v_2\cdots$ with $|v_k|=k$. If $A=\{\textsl{o},v_1,v_2,\ldots\}$
then the characteristic function $\chi_{A}$, namely, $\chi_{A}(v)=1$ for $v\in A$ and zero elsewhere, belongs to
$\mathbb{T}_{p}$ for all $0<p<\infty$ but is not in $\mathcal{L}_w$. This concludes the proof that $\mathcal{L}_w$ is not
comparable with $\mathbb{T}_{p}$ for all $0<p<\infty$.

\item
Clearly, $\mathbb{T}_{\infty} \subseteq \big (\bigcap\limits_{0<p<\infty}\mathbb{T}_{p}\big ) \bigcap \mathcal{L}$,
whereas $\mathcal{L}_w$ is not comparable with $\mathbb{T}_{\infty}$. For 2-homogeneous trees, this inclusion relation
becomes an equality. This is because of the fact that there is no unbounded function in $\mathbb{T}_{p}$ for $2$-homogeneous
trees, which can be observed from the definition of $\mathbb{T}_{p}$.
\end{enumerate}
\erems

\subsection{Separability}
In order to state results about the separability of $\mathbb{T}_{p,0}$ and $\mathbb{T}_{p}$, we need to introduce few notations. Denote by
$C_{c}(T)$ the set of all functions $f\colon T \to\mathbb{C}$ such that $M_{p}(n,f)=0$ for all but finitely many $n'$s.
Also the closure of $C_{c}(T)$ under $\|.\|_p$ is denoted by $\overline{C_{c}(T)}$.

\blem
For $0<p\leq\infty$, we have $\overline{C_{c}(T)}= \mathbb{T}_{p,0}$.
\elem
\bpf
Let $f\in \mathbb{T}_{p,0}$ and, for each $n$, define $\{f_n\}$ by
$$ f_n(v)=\left\{
\begin{array}{rl}
f(v) & \mbox{ if } |v|\leq n\\
0 & \mbox{ otherwise}.
\end{array}
\right.
$$
Clearly, $f_n \in C_{c}(T)$ for each $n\in \IN$ and
$$M_p(k,f-f_n)=\left\{
\begin{array}{rl}
M_p(k,f) & \mbox{ if } k>n \\
0 & \mbox{ otherwise.}
\end{array}
\right.
$$
Therefore, we see that
$$\|f-f_n\|_p = \sup\limits_{m\in \mathbb{N}_{0}} M_{p}(m,f-f_n)= \sup\limits_{m>n} M_{p}(m,f)
$$
and, because $M_{p}(m,f)\ra0$ as $m\ra\infty$, it follows that $\|f-f_n\|_p\ra0$ as $n\ra\infty$.
This completes the proof.
\epf

\bthm\label{thm:sepp0}
For $0<p\leq\infty$, $\mathbb{T}_{p,0}$ is a separable space.
\ethm
\bpf
It is easy to verify that $B=\{\chi_{\{v\}} :\, v\in T\}$ is a basis for $C_{c}(T)$.
Since $B$ is countable, $C_{c}(T)$ is separable. Since $C_{c}(T)$ is dense in $\mathbb{T}_{p,0}$, we conclude that $\mathbb{T}_{p,0}$ is separable
and the theorem follows.
\epf

We remark that $C_c(T)$ cannot be a Banach space with respect to any norm, since it has a countably infinite basis.

\bthm\label{thm:separable}
For $0<p\leq\infty$, $\mathbb{T}_{p}$ is not separable.
\ethm
\bpf
Let $E\subset \mathbb{T}_{p} $ denote the set of all radial constant functions $f$ whose range is subset of $\{0,1\}$. Let $f,g \in E$ and$f\neq g$.
Then there exists a $v\in T$ such that $f(v)\neq g(v)$. Since $f,g \in E$, we have $M_p(n,f-g)\leq 1$ for all $n$. On the other hand, $M_p(|v|,f-g)=1$
and hence, $\|f-g\|_p= \sup M_p(n,f-g)=1$. It is easy to check that $E$ is an uncountable subset of $\mathbb{T}_{p}$.
Since any two distinct elements of $E$ must be of distance $1$ apart and $E$ is uncountable, it follows that any dense subset
of $\mathbb{T}_{p}$ cannot be countable. Consequently, $\mathbb{T}_{p}$ is not a separable space.
\epf

\subsection{Growth estimate and consequences}

\blem \label{lem:bound}
Let $T$ be a $q+1$ homogeneous tree rooted at $\textsl{o}$ and $0<p<\infty$. Then, for $v\in T \setminus \{\textsl{o}\}$, we have the following:
\begin{enumerate}
\item[(a)] If $f\in \mathbb{T}_{p}$, then
$|f(v)|\leq \{(q+1)q^{|v|-1}\}^{\frac{1}{p}} \|f\|_p. $
\item[(b)] If $f\in \mathbb{T}_{p,0}$, then
$$ \lim_{|v| \ra \infty} \frac{f(v)}{\{(q+1)q^{|v|-1}\}^{\frac{1}{p}}} =0.
$$
\end{enumerate}
The results are sharp.
\elem
\bpf Fix $v\in T \setminus \{\textsl{o}\}$ and let $n=|v|$. Then,
$$|f(v)|^p\leq \sum\limits_{|w|=n}|f(w)|^{p}= (q+1)q^{n-1} M_p^p(n,f)
$$
so that
$|f(v)|\leq \{(q+1)q^{n-1}\}^{\frac{1}{p}} M_p(n,f)
$
and thus,
$$ \frac{|f(v)|}{\{(q+1)q^{n-1}\}^{\frac{1}{p}}} \leq M_p(n,f) \leq \|f\|_p.
$$
The desired results follow.

In order to prove the sharpness, we fix $v\in T \setminus \{\textsl{o}\}$. Define $f(v)=\{(q+1)q^{|v|-1}\}^{\frac{1}{p}}$ and
$0$ elsewhere. We now let $m=|v|$ so that $M_p(n,f)=0$ for every $n\neq m$ and
$$M_p(m,f)= \left (\frac{1}{(q+1)q^{m-1}}\sum\limits_{|w|=m}|f(w)|^{p} \right )^{\frac{1}{p}}=
\left (\frac{1}{(q+1)q^{m-1}}|f(v)|^{p} \right )^{\frac{1}{p}}=1.
$$
We obtain that $\|f\|_{p} = \sup\limits_{n\in \mathbb{N}_{0}} M_{p}(n,f)=1$ and hence,
$$|f(v)|=f(v)=\{(q+1)q^{|v|-1}\}^{\frac{1}{p}} \|f\|_p.
$$
We conclude the proof.
\epf

Lemma \ref{lem:bound}(a) clearly holds if $v=\textsl{o}$, because $|f(\textsl{o})| \leq \|f\|_p.$

\bprop\label{prop1}
(Compare with \cite[Corollary, p.~10]{Shapiro:Book})
Convergence in $\|.\|_p$ $(0<p\leq\infty)$ implies uniform convergence  on compact subsets of $T$.
\eprop
\bpf
The edge counting distance on $T$ induces the discrete metric. So, finite subsets are the only compact sets in $T$.
Let $K$ be an arbitrary compact subset of $T$. Then there exists an $N\in \mathbb{N}$ such that $|v|\leq N$ for every $v\in K$.
The proposition trivially holds for the case $p=\infty$, because given a function $f$ and a sequence $\{f_n\}$ converging to $f$ in norm,
$$\sup_{v\in K}|(f_n-f)(v)|\leq \|f_n-f\|_\infty.
$$
Next, we consider the case $0<p<\infty$. From Lemma \ref{lem:bound}, given a function $f$ in $\mathbb{T}_{p}$, we have
$$|f(v)|\leq \{(q+1)q^{|v|-1}\}^{\frac{1}{p}} \|f\|_p \leq \{(q+1)q^{N-1}\}^{\frac{1}{p}} \|f\|_p \mbox{ for every $v\in K$}.
$$
This gives
$$\sup_{v\in K}|f(v)|\leq \{(q+1)q^{N-1}\}^{\frac{1}{p}} \|f\|_p
$$
and thus, by replacing $f$ by $f_n-f$, we conclude that convergence in $\|.\|_p$ implies uniform convergence  on compact subsets of $T$.
\epf

 Uniform convergence on compact subsets of $T$ does not necessarily imply the convergence in $\|.\|_p ~(0<p\leq\infty)$ as can be seen from
the following example.

\beg
Consider the function $f\equiv 1$. For each $n$, define $\{f_n\}$ by
$$f_n(v):=\left\{
\begin{array}{rl}
f(v)\,(=1) & \mbox{ if } |v|\leq n\\
0 & \mbox{ otherwise.}
\end{array}
\right.
$$
Then
$$M_p(k,f-f_n)=\left\{
\begin{array}{rl}
M_p(k,f)\, (=1) & \mbox{if } k>n\\
0 & \mbox{ otherwise}.
\end{array}
\right.
$$
Let $K$ be a compact subset of $T$. Then there exists an $N\in \mathbb{N}$ such that $|v|\leq N$ for every $v\in K$ and
$\sup_{v\in K}|(f_n-f)(v)|=0$ for every $n>N$. It follows that $\{f_n\}$ converges uniformly on compact subsets of $T$ to $f$.
On the other hand,
$$\|f-f_n\|_p = \sup\limits_{m\in \mathbb{N}_{0}} M_{p}(m,f-f_n)= 1 ~\mbox{ for every $n\in \IN$.}
$$
Hence, $\{f_n\}$ does not converge to $f$ in $\|.\|_p$.
\eeg

From Proposition \ref{prop1} and the above remark, we observe that the the topology of uniform convergence on the compact
subsets of $T$ on $\mathbb{T}_{p}$ is similar to that of analytic cases such as $H^p$ spaces. This observation raises a natural question.
Is $\mathbb{T}_{p}$ complete in the topology of uniform convergence on compact sets?
The following example shows that the answer is negative. For each $n$, define $\{f_n\}$ by
$$f_n(v):=\left\{
\begin{array}{rl}
|v| & \mbox{ if } |v|\leq n\\
0 & \mbox{ otherwise.}
\end{array}
\right.
$$
Let $K$ be a compact subset of $T$. Then there exists an $N\in \mathbb{N}$ such that $|v|\leq N$ for every $v\in K$.
For $N<n<m$, $f_n(v)=f_m(v)$ for all $v\in K$. It is easy to see that $\{f_n\}$ is a Cauchy sequence in the topology of uniform
convergence on compact sets and $\{f_n\}$ converges pointwise to the function $f(v)=|v|$.
Note that $f$ can be the only possible limit of $\{f_n\}$ in the topology of uniform convergence on compact sets.
Since $\mathbb{T}_{p}$ contains the sequence $\{f_n\}$ but not $f$, $\mathbb{T}_{p}$ cannot be complete under the topology of
uniform convergence on compact sets.

\section{Multiplication Operators on $\mathbb{T}_{p}$ and $\mathbb{T}_{p,0}$} \label{Multiplication Operators}

We now recall the following definitions.

Let $X$ be a complex normed linear space consisting of complex valued functions defined on a set $\Omega$. If $\psi$ is
a complex valued function defined on $\Omega$, then the multiplication operator with symbol $\psi$ is defined by $M_\psi f = \psi f$ for every $f \in X$.

A Banach space $X$ on $\Omega$ said to be a \textsl{functional Banach space} if for each $v\in \Omega$, the point
evaluation map $e_v : f\in X \mapsto f(v)$ is a bounded linear functional on $X$. The following result is well-known.

\blem \cite[Lemma 11]{Duren:functinal} \label{lem:functional}
Let $X$ be a functional Banach space on the set $\Omega$ and $\psi$ be a complex valued function on $\Omega$ such that $M_\psi$ maps
$X$ into itself. Then $M_\psi$ is bounded on $X$ and $|\psi(v)|\leq \|M_\psi\|$ for all $v\in \Omega$. In particular, $\psi$ is bounded function.
\elem

It is natural to ask whether $\mathbb{T}_{p}$ and $\mathbb{T}_{p,0}$ are functional Banach spaces.

\subsection{Boundedness}

\bprop
For $1\leq p\leq\infty$, $\mathbb{T}_{p}$ and $\mathbb{T}_{p,0}$ are functional Banach spaces.
\eprop
\bpf
First we consider the case when $p=\infty$. Since $|e_v(f)|=|f(v)|\leq \|f\|_\infty$ for every $v \in T$, it follows that
$\mathbb{T}_{\infty}$ is a functional Banach space.

Let us now consider the case when $1\leq p<\infty$. The point evaluation map $e_\textsl{o}$ is a bounded linear functional on $\mathbb{T}_{p}$
because of the fact that $|f(\textsl{o})|\leq \|f\|_p$ for every $f$ in $\mathbb{T}_{p}$.
Now, we fix $v\in T$ and $v\neq\textsl{o}$. Then, from Lemma \ref{lem:bound}, we have
$$|e_v(f)|=|f(v)|\leq \{(q+1)q^{|v|-1}\}^{\frac{1}{p}} \|f\|_p ~\mbox{ for every }~ f\in \mathbb{T}_{p}.
$$
So $e_v$ is a bounded linear functional on $\mathbb{T}_{p}$ with $\|e_v\|\leq \{(q+1)q^{|v|-1}\}^{\frac{1}{p}}$. Hence $\mathbb{T}_{p}$ is a
functional Banach space. A similar proof works also for the space $\mathbb{T}_{p,0}$. The proof is complete.
\epf

\bthm \label{thm:boundedness}
Let $T$ be a  $q+1$ homogeneous tree rooted at $\textsl{o}$, $\psi$ be a function on $T$.
Let $X$ be $\mathbb{T}_{p}$ or $\mathbb{T}_{p,0}$ or $C_c(T)$ with norm $\|.\|_p$, where $0<p\leq\infty$.
Then the following are equivalent (compare with \cite [Proposition~2]{Dragan:survey}).
\begin{enumerate}
\item[(a)] $M_\psi$ is a bounded linear operator from $X$ to $X$,
\item[(b)] $\psi$ is a bounded function on $T$, i.e., $\psi \in \mathbb{T}_{\infty}$.
\end{enumerate}
\ethm
\bpf
Let $X$ be $\mathbb{T}_{p}$ or $\mathbb{T}_{p,0}$ or $C_c(T)$ with norm $\|.\|_p$, where $0<p\leq\infty$.\\

${\rm (a)}\Rightarrow {\rm (b)}:$ We will prove this implication by contradiction. Suppose that $\psi$ is an unbounded function on $T$.
Then there exists a sequence of vertices $\{v_k\}$ such that
$$|v_1|<|v_2|<|v_3|\cdots ~\mbox{ and }~ |\psi(v_k)|\geq k.
$$
For each $k$, define $f_k:T\rightarrow \mathbb{C}$ by $f_k= C_{k,p} \chi_{\{v_k\}}$, where the constants $C_{k,p}$'s are
chosen in such a way that $\|f_k\|_p=1$. Note that $f_k \in X$ for each $k \in \mathbb{N}$. We obtain that
$$k=k \|f_k\|_p \leq |\psi(v_k)| \|f_k\|_p = M_p(|v_k|,\psi f_k) \leq \|\psi f_k\|_p = \|M_\psi f_k\|_p
$$
for every $k \in \IN$, which gives a contradiction to our assumption. Hence, $\psi$ is a bounded function on $T$.\\

${\rm (b)}\Rightarrow {\rm (a)}:$ Suppose that $\psi$ is a bounded function on $T$ and $0< p<\infty$. Then for any $f\in X$,
\beqq
M_p(n,M_\psi f)&= & \left (\frac{1}{(q+1)q^{n-1}}\sum\limits_{|v|=n}|(M_\psi f)(v)|^{p}\right)^{\frac{1}{p}}\\
&= & \left (\frac{1}{(q+1)q^{n-1}}\sum\limits_{|v|=n}|\psi(v)|^p |f(v)|^{p} \right )^{\frac{1}{p}},
\eeqq
which shows that
\be \label{eq:integ-mean}
M_p(n,M_\psi f) \leq \|\psi\|_\infty M_p(n,f).
\ee
For $p=\infty$, the inequality (\ref{eq:integ-mean}) is trivially holds. From (\ref{eq:integ-mean}), one can also observe that
$M_\psi f \in X$ whenever $f \in X$. Taking the supremum over $n \in \mathbb{N}_{0}$ on both sides of (\ref{eq:integ-mean}),
we deduce that $\| M_\psi f \|_{p} \leq \|\psi\|_\infty \|f\|_{p}$ and thus, $M_\psi$ is bounded linear operator from $X$ to $X$
with $\|M_\psi\| \leq \|\psi\|_\infty $.
\epf

\brem
Let $X$ be $\mathbb{T}_{p}$ or $\mathbb{T}_{p,0}$ with norm $\|.\|_p$, where $1\leq p\leq\infty$. In this case, $X$ becomes a functional
Banach space. Then by Lemma \ref{lem:functional}, one has $|\psi(v)|\leq \|M_\psi\|$ for all $v\in T$ which by taking the
supremum gives $\|\psi\|_\infty  \leq \|M_\psi\|$. Therefore, by Theorem \ref{thm:boundedness}, it follows that
$\|M_\psi\| = \|\psi\|_\infty $.
\erem

\subsection{Spectrum}

Let $X$ be a normed linear space and $A$ be a bounded linear operator on $X$. The \textit{point spectrum} $\sigma_e(A)$ of $A$
consists of all $\lambda \in \IC$ such that $A-\lambda I$ is not injective. Thus $\lambda \in \sigma_e(A)$ if and only if there is
some nonzero $x\in X$ such that $A(x)=\lambda x$. The \textit{approximate point spectrum} $\sigma_a(A)$ of $A$ consists of all
$\lambda \in \IC$ such that $A-\lambda I$ is not bounded below. Thus $\lambda \in \sigma_a(A)$ if and only if there is a sequence
$\{x_n\}$ in $X$ such that $\|x_n\|=1$ for each $n$ and $\|A(x_n)-\lambda x_n\| \ra 0$ as $n \ra \infty$.
The \textit{spectrum} $\sigma(A)$ of $A$ consists of all $\lambda \in \IC$ such that $A-\lambda I$ is not invertible.

It is clear from the definition that
\be
\sigma_e(A)\subseteq \sigma_a(A) \subseteq \sigma(A).
\ee

It is well-known that the spectrum of a bounded linear operator $A$ on a Banach space $X$ over $\IC$ is a nonempty compact
subset of $\IC$ (see \cite[Theorem 3.6]{Conway:Book}). Also, the boundary of the spectrum is contained in the approximate
point spectrum of $A$ (\cite[Proposition 6.7]{Conway:Book}) and approximate point spectrum is a closed subset of $\IC$.

\bthm
For $1\leq p \leq \infty$, let $X$ equal either $\mathbb{T}_{p}$ or $\mathbb{T}_{p,0}$ or $C_c(T)$, and let
$M_\psi$ be a bounded multiplication operator on $X$ with norm $\|.\|_p$. Then
\begin{enumerate}
\item[(a)] $\sigma_e(M_\psi)=$ Range of $\psi = \psi(T)$;
\item[(b)] $\sigma(M_\psi)= \sigma_a(M_\psi)= \overline{\psi(T)}$.
\end{enumerate}
\ethm
\bpf
In order to prove (a), we begin by letting $\lambda \in \sigma_e(M_\psi)$. Then there exists a nonzero function $f\in X$ such that
$\psi f = M_\psi f = \lambda f$. Since $f\neq0$, there is a vertex $v$ such that $f(v)\neq0$ and $(\psi(v)-\lambda)f(v)=0$.
Thus, $\lambda = \psi(v) \in \psi(T)$ and therefore, $\sigma_e(M_\psi)\subseteq \psi(T)$.

Conversely, suppose that $\alpha \in \psi(T)$. Then, there exists a vertex $v$ such that $\psi(v)= \alpha$. Thus,
$M_\psi( \chi_{\{v\}} )= \alpha \chi_{\{v\}}$ and $0 \neq \chi_{\{v\}} \in X$. This gives $\alpha \in \sigma_e(M_\psi)$.
Hence, $\sigma_e(M_\psi)= \psi(T)$.

Before proving (b), we observe that for every $\lambda \in \IC, M_\psi-\lambda I = M_{\psi-\lambda}$. Thus, $\lambda \in \sigma(M_\psi)$
if and only if $M_{\psi-\lambda}$ is not invertible if and only if $\frac{1}{\psi-\lambda}$ is not a bounded function on $T$ (see Theorem \ref{thm:boundedness}).
Now, we let $\lambda \notin \overline{\psi(T)}$. Since the complement of $\overline{\psi(T)}$ is open,
there exists an $r>0$ such that disk of radius $r$ centered at $\lambda$ is a subset of $\mathbb{C} \setminus \overline{\psi(T)}$.
So, $|\psi(v)-\lambda| \geq r$ for every $v\in T$. Thus $\frac{1}{\psi-\lambda}$ is a bounded function and therefore,
by Theorem \ref{thm:boundedness}, $M_{\frac{1}{\psi-\lambda}}$ is a bounded operator on $X$. It is easy to verify that
$\ds M_{\frac{1}{\psi-\lambda}}$ is the inverse of $\ds M_{\psi-\lambda}$ and hence, $\ds M_{\psi-\lambda}$ is invertible.
We conclude that $\lambda$ cannot be in the spectrum, which in turn implies that
$\sigma(M_\psi) \subseteq \overline{\psi(T)}$.

On the other hand, $\psi(T) = \sigma_e(M_\psi)\subseteq \sigma_a(M_\psi)\subseteq\sigma(M_\psi) \subseteq \overline{\psi(T)}$ and the fact that the
approximate point spectrum and spectrum are closed subsets of $\IC$ give that $\sigma(M_\psi)= \sigma_a(M_\psi)= \overline{\psi(T)}$.
\epf

\brem
Let $X$ be $\mathbb{T}_{p}$ or $\mathbb{T}_{p,0}$ or $C_c(T)$ with $\|.\|_p$, where $1\leq p\leq\infty$. Then
$M_\psi: X \rightarrow X$ is not injective if and only if $0 \in \sigma_e(M_\psi) = \psi(T)$. So, $0$ is in the range of $\psi$
is a necessary and sufficient condition for $M_\psi$ not being injective on $X$.
\erem

\subsection{Compactness}

\bthm \label{thm:Compactness}
Let $X$ be either $\mathbb{T}_{p}$ or $\mathbb{T}_{p,0}$, where $1\leq p\leq\infty$ and let $M_\psi: X \rightarrow X$ be a
bounded multiplication operator on $X$. Then $M_\psi$ is a compact operator on $X$ if and only if $\psi(v)\rightarrow 0$ as $|v|\rightarrow \infty$.
\ethm
\bpf
Let $M_\psi$ be a compact operator on $X$. Then, from \cite[Theorem 4.25]{Rudin:Book}, $\sigma_e(M_\psi)=\psi(T)$ (as well as $\sigma(M_\psi)$)
is a countable set with $0$ as the only possible limit point. Suppose $\psi(v)\not\rightarrow 0$ as $|v|\rightarrow \infty$.
Then there exists an $\epsilon>0$ and a sequence $\{v_k\}$ in $T$ such that $|v_k|\rightarrow \infty$ and $|\psi(v_k)| \geq \epsilon$ for all $k$.
By Bolzano-Weierstrass theorem, $\{\psi(v_k)\}$ has a limit point. Because $|\psi(v_k)| \geq \epsilon$ for all $k$, we have a contradiction to the fact that
$0$ is the only possible limit point. Hence, $\psi(v)\rightarrow 0$ as $|v|\rightarrow \infty$.

For the proof of the converse part, we use the fact that the set of compact operators on $X$ is a closed subspace of the set of all bounded operators on $X$
(see \cite[Theorem 4.18 part (c)]{Rudin:Book}).

Consider a function $\psi$ from $C_c(T)$. Then there exists an $N\in \mathbb{N}$ such that
$\psi(v)=0$ for every $|v|>N$. So, $(M_\psi f)(v)=\psi(v) f(v)=0$ for every $|v|>N$ and for every $f\in X$. Thus, the range of
$M_\psi$ is a finite dimensional subspace, which shows that $M_\psi$ is a compact operator (\cite[Theorem 4.18 part\,(a)]{Rudin:Book}).

Let $\psi$ be a arbitrary function such that $\psi(v)\rightarrow 0$ as $|v|\rightarrow \infty$.
For each $n$, define $\{f_n\}$ by
$$ f_n(v):=\left\{
\begin{array}{rl}
\psi(v) & \mbox{ if } |v|\leq n\\
0 & \mbox{ otherwise}.
\end{array}
\right.
$$
By definition $f_n \in C_c(T)$ for every $n$. Thus, $M_{f_n}$ is a compact operator for every $n$. Moreover,
$$\|M_{f_n}-M_\psi\|= \|M_{f_n-\psi}\|= \|f_n-\psi\|_\infty= \sup\limits_{|v|>n} |\psi(v)|,
$$
which approaches to zero as $|v|\rightarrow \infty$, because $\psi(v)\rightarrow 0$ as $|v|\rightarrow \infty$.
Thus, $M_\psi$ is the limit (in the operator norm) of a sequence $\{M_{f_n}\}$ of compact operators, and hence,
$M_\psi$ is compact on $X$. The proof is now complete.
\epf

\blem \label{lem:Compactness}
Let $X$ be either $\mathbb{T}_{p}$ or $\mathbb{T}_{p,0}$, where $1\leq p\leq\infty$ and let $M_\psi: X \rightarrow X$ be a bounded multiplication operator on $X$.
If $M_\psi$ is a compact operator on $X$, then, for every bounded sequence $\{f_n\}$ in $X$ converging to $0$ pointwise, the sequence
$\|\psi f_n\|\rightarrow 0$ as $n\rightarrow\infty$.
\elem
\bpf
Suppose that $\{g_n\}$ in $X$ is a bounded sequence converging to $0$ pointwise. Since $M_\psi$ is a compact operator, there is a subsequence
$\{g_{n_k}\}$ of $\{g_{n}\}$ such that $\{\psi g_{n_k}\} = \{M_\psi(g_{n_k})\}$ converges in $\|.\|_p$ to some function, say, $g$.
It follows that $\{\psi g_{n_k}\}$ converges to $g$ pointwise. Since the convergence of $\{g_n\}$ to $0$ implies that $g \equiv 0$, we deduce that
$\{\psi g_{n_k}\}$ converges to $0$ in $\|.\|_p$.

Let $\{f_n\}$ be a bounded sequence in $X$ converging to $0$ pointwise. We claim that $\|M_\psi(f_n)\| = \|\psi f_n\|\rightarrow 0$
as $n\rightarrow\infty$. Suppose that $\|M_\psi(f_n)\| \not\rightarrow 0$ as $n \rightarrow \infty$.
Then there exists a subsequence $\{f_{n_j}\}$ and an $\epsilon >0$ such that $\|M_\psi(f_{n_j})\| \geq \epsilon$ for all $j$.
By taking $g_{n} = f_{n_j}$ in the last paragraph, we find that $\{\psi g_{n_k}\}$ converges to $0$ in $\|.\|_p$,
which is not possible because $\|M_\psi(f_{n_j})\| \geq \epsilon$ for all $j$. Hence, $\|M_\psi(f_n)\| = \|\psi f_n\|\rightarrow 0$
as $n\rightarrow\infty$, and the proof is complete.
\epf

It is well-known (cf. \cite[Proposition 3.3]{Conway:Book}) that if $A$ is a compact operator on a normed linear space $Y$, then
the image of every weak convergent sequence is norm convergent,
i.e., $x_n \rightarrow x$ weakly implies that $A x_n \rightarrow Ax.$ In view of this result and Lemma \ref{lem:Compactness},
we have a natural question:
\textit{Is there a relationship between weak convergence to $0$ in $X$ and the boundedness of a sequence in $X$ converging to $0$ pointwise?}

\subsection{Upper bound for the essential norm}

Let $CL(X)$ denote the set of all compact operators on $X$. The essential norm $\|A\|_e$ of an bounded operator $A$ on $X$ defined
to be the distance between $A$ and $CL(X)$:
$$\|A\|_e=\dist(A,CL(X))=\inf\{\|A-K\|: \, K\in CL(X)\}.
$$
The following theorem is a natural generalization of Theorem \ref{thm:Compactness}.

\bthm \label{thm-upessnorm}
Let $M_\psi$ be a bounded multiplication operator on $\mathbb{T}_{p}$ or on $\mathbb{T}_{p,0}$, where $1\leq p\leq\infty$. Then
$$\|M_\psi\|_e \leq \limsup\limits_{n\rightarrow\infty} M_\infty(n,\psi)=\lim\limits_{n\rightarrow\infty} \sup\limits_{|v|\geq n} |\psi(v)|.
$$
\ethm
\bpf
For each $n\in \mathbb{N}$, define $\{\psi_n\}$ by
$$ \psi_n(v):=\left\{
\begin{array}{rl}
\psi(v) & \mbox{ if } |v|<n\\
0 & \mbox{ otherwise}.
\end{array}
\right.
$$
Clearly, $M_{\psi_n}$ is a compact operator for every $n\in \mathbb{N}$ and thus, for every $n\in \mathbb{N}$,
\beqq
\|M_\psi\|_e &=& \inf\{\|M_\psi-K\|:K\in CL(X)\} \\
&\leq & \|M_\psi-M_{\psi_n}\| = \|\psi-\psi_n\|_\infty \\
&=& \sup\limits_{|v|\geq n} |\psi(v)| = \sup\limits_{m\geq n} M_\infty(m,\psi).
\eeqq
Hence, $\|M_\psi\|_e \leq \inf\limits_{n\in \mathbb{N}} \sup\limits_{|v|\geq n} |\psi(v)| = \limsup\limits_{n\rightarrow\infty} M_\infty(n,\psi)$
and the proof is complete.
\epf

We would like to point out that one way implication of Theorem \ref{thm:Compactness} follows from Theorem \ref{thm-upessnorm}. Indeed if
if $\psi(v)\rightarrow 0$ as $|v|\rightarrow \infty$, then Theorem \ref{thm-upessnorm} gives that $\|M_\psi\|_e =0$ and thus,
$M_\psi$ is a compact operator.

\subsection{Isometry}

A map $A$ from a normed linear space $(X,\|.\|_X)$ to a normed linear space $(Y,\|.\|_Y)$ is said to be an isometry
if $\|Ax\|_Y = \|x\|_X$ for all $x\in X$.

\bthm
Let $X$ be either $\mathbb{T}_{p}$ or $\mathbb{T}_{p,0}$, where $0<p\leq\infty$ and let $M_\psi: X \rightarrow X$ be a bounded multiplication
operator on $X$. Then $M_\psi$ is an isometry on $X$ if and only if $|\psi(v)|=1$ for all $v\in T$.
\ethm
\bpf
Suppose that $|\psi(v)|=1$ for all $v\in T$. Then $M_p(n,\psi f)=M_p(n, f)$ for all $n$, which shows that $\|f\|_p=\|\psi f\|_p=\|M_\psi(f)\|_p$ and
thus, $M_\psi$ is an isometry on $X$.

Conversely, suppose that $M_\psi$ is an isometry on $X$. First we consider the case $p=\infty$. Let $f$ be $\chi_{\{v\}}$. Because
$M_\psi$ is an isometry on $X$, we have
$$|\psi(v)|=\|M_\psi(f)\|_\infty=\|\psi f\|_\infty=\|f\|_\infty=1,
$$
which holds for every $v\in T$. Hence $|\psi(v)|=1$ for all $v\in T$.

Next, we consider the case $0<p<\infty$. Let $f$ be $\chi_{\{\textsl{o}\}}$. Since $M_\psi$ is an isometry, $|\psi(\textsl{o})|=1$.
Take an arbitrary element $v\in T$ with $|v|\geq 1$ and let $f$ be $\chi_{\{v\}}$. Moreover, since $M_\psi$ is an isometry, we have
$$ \left (\frac{|\psi(v)|^{p}}{(q+1)q^{n-1}} \right )^{\frac{1}{p}}=\|M_\psi(f)\|_p=\|f\|_p=\left (\frac{1}{(q+1)q^{n-1}} \right )^{\frac{1}{p}},
$$
which shows that $|\psi(v)|=1$ for all $v\in T$.
\epf

\subsection*{Acknowledgement}
The authors thank the referee for many useful comments which improve the presentation considerably.
The first author thanks the Council of Scientific and Industrial Research (CSIR), India,
for providing financial support in the form of a SPM Fellowship to carry out this research.
The second  author is on leave from the Indian Institute of Technology Madras, India.

\end{document}